\documentclass[reqno,11pt]{amsart}
\numberwithin{equation}{section}
\usepackage{amsmath}
\usepackage{esint} % permette di fare l'integrale tagliato
\usepackage{amsthm}
\usepackage{epsfig}
\usepackage{psfrag}
\usepackage{graphicx}
\usepackage{graphpap,latexsym,epsf}
\usepackage{color}
\usepackage{amssymb,mathrsfs,enumerate}
\usepackage[colorlinks, citecolor=citegreen, linkcolor=refred]
{hyperref}
\newcommand{\R}{\mathbb{R}}
\newcommand{\RRR}{\mathbb{R}^{3{\times}3}}

\newcommand{\Sph}{\mathbb{S}}

\newcommand{\pa}{\partial}
\newcommand{\na}{\nabla}
\newcommand{\Om}{\Omega}

\newcommand{\ep}{\varepsilon}

\newcommand{\dd}{{\rm d}}
\newcommand{\SO}{{\rm SO}}

\newcommand{\tr}{{\rm tr}}
\newcommand{\sym}{{\rm sym}\,}

\newcommand{\kk}{k}

\mathchardef\emptyset="001F

\definecolor{vgreen}{rgb}{0.1,0.5,0.2}
\definecolor{viola}{RGB}{85,26,139}
\definecolor{citegreen}{rgb}{0,0.6,0}
\definecolor{refred}{rgb}{0.8,0,0}
\newtheorem{theorem}{Theorem}[section]
\newtheorem{remark}[theorem]{Remark}

\begin{document}
\title[Dimension reduction via $\Gamma$-convergence 
for soft active materials]
{Dimension reduction via $\Gamma$-convergence\\ 
for soft active materials}

\author[V.~Agostiniani]{Virginia Agostiniani}
\address{SISSA, via Bonomea 265, 34136 Trieste - Italy} 
\email{vagostin@sissa.it}

\author[A.~DeSimone]{Antonio DeSimone}
\address{SISSA, via Bonomea 265, 34136 Trieste - Italy}
\email{desimone@sissa.it}

\begin{abstract} 
We present a rigorous derivation of
dimensionally reduced theories for thin sheets of nematic elastomers, in the finite bending regime.
Focusing on the case of \emph{twist} nematic texture, we obtain 
2D and 1D models for wide and narrow ribbons exhibiting spontaneous flexure and torsion.
We also discuss some variants to the case of twist nematic texture, 
which lead to 2D models with different \emph{target curvature} tensors.
In particular, we analyse cases where the nematic texture
leads to zero or positive
Gaussian target curvature,
and the case of bilayers.
\end{abstract}

\maketitle

%\date{\today}

% \noindent\textsc{MSC (2010): 74B20 (nonlinear elasticity),
% 76A15 (liquid crystals), 74K10 (rods), 74K20 (plates)}
% \smallskip
% \noindent\keywords{nematic elastomers, dimension reduction, 
% plate theory, rod theory} 

%%%%%%%%%%%%%%%%%%%%%%%%%%%%%%%%%
%%%%%%%%%%%%%%%%%%%%%%%%%%%%%%%%%

%%%%%%%%%%%%%%%%%%%%%%%%%%%%%%%%%%%%%%%%%%

\section{Introduction}

%%%%%%%%%%%%%%%%%%%%%%%%%%%%%%%%%%%%%%%%%%

We discuss in this paper some recent (and some new) results
motivated by shape programming problems for thin structures,
made of soft active materials.
In particular, we focus our attention on nematic elastomers
\cite{Ag_Bl_Ko,AgDalDe2,AgDe2,ferro,DeSim_Dolz,Warner_book}.
These are polymeric materials that, thanks to the coupling
between elasticity and orientational nematic order, 
undergo spontaneous deformations as a consequence
of a temperature-driven phase transformation.

When a nematic texture (i.e.\ a spatially-dependent
orientation of the nematic director) is present,
non-uniform spontaneous strains lead to stress build-up
and non-trivial configurational changes.
Here, we explore in particular the consequences of
through-the-thickness changes of orientation in thin sheets.
Plates with spontaneous and controllable 
\emph{curvature} emerge in 
this way. When the plate is in the form of a narrow
ribbon, 
rods with spontaneous \emph{flexure} and 
\emph{torsion} are generated.
Similar issues have been the object of recent intense investigation,
see e.g. 
\cite{Ag_De_bend,Ago_DeS_Kou,hybrid,Sawa2011,Urayama2013}
and also
\cite{Aha_Sha_Kup,Bhat_199,Bla_Ter_War,CDD2,CCDD1,
Lucantonio2016,Teresi_Varano}.
The point of view we adopt here is that of a systematic 
derivation of dimensionally-reduced models
by means of $\Gamma$-convergence.
Far from being just a technical exercise,
this technique allows one to \emph{derive},
rather than \emph{postulate},  
the functional form and the material parameters 
(elastic constants, spontaneous curvature, etc.) 
for dimensionally reduced models of these structures.

Building upon our recent work \cite{Ag_De_bend,Ago_DeS_Kou},
which has been in turn inspired by 
\cite{Freddi2015,F_J_M_2002,Schmidt2007},
we discuss in this paper the derivation of models capable of describing 
non-trivial shapes that can be spontaneously exhibited
by thin sheets of nematic elastomers.
In particular, we focus on the \emph{twist} director geometry,
and we present a derivation of reduced models for wide and narrow ribbons
(namely, rods whose cross-section is a thin rectangle)
by following a sequential dimension reduction technique, namely, a 3D$\rightarrow$2D$\rightarrow$1D procedure, 
see Section \ref{sec:main}.
Plates made of twist nematic elastomers spontaneously
deform into cylindrical ribbons. 
%(i.e., ribbons \emph{wrapped} around cylinders). 
This geometry emerges as a ``compromise'' between the
energetic drive towards a non-zero \emph{target curvature}
(with negative Gaussian curvature) 
and an isometry constraint on the mid-surface
(which rules out surfaces with non-zero Gaussian curvature).
Narrow ribbons obtained from rectangular plates of 
vanishing small width behave as rods with spontaneous 
\emph{flexure} and \emph{torsion}.

Moreover, we also consider in Section  \ref{sec:var} some alternative scenarios to the case of twist nematic texture, in which the target curvature tensor 
of the 2D reduced model has either zero or positive Gaussian curvature. 
Nematic textures realising these possibilities are the splay-bend one and uniform director alignment perpendicular to the mid-plane, respectively.
Alternatively, the three alternative scenarios (positive, negative, or zero Gaussian curvature of the target curvature tensor in the 2D plate theory) can be realised in bilayers, by tuning the spontaneous strains in the two halves of the bilayer.
In all these cases, configurations with non-zero Gaussian curvature are again frustrated by the isometry constraint.
Hence, the configurations that are spontaneously exhibited 
are always ribbons wrapped around cylinders.
Much of the material relating to nematic elastomer sheets draws on our analyses contained in \cite{Ag_De_bend,Ago_DeS_Kou},  but it is presented here from the unifying perspective of a sequential dimensional reduction, first from three to two dimensions, and then from two to one dimension. The analysis of thin bilayer sheets is new.
We refer the reader to \cite{Nardinocchi} 
for a discussion of bilayer beams.

%%%%%%%%%%%%%%%%%%%%%%%%%%%%%%%%%%%%%%%%%%
%%%%%%%%%%%%%%%%%%%%%%%%%%%%%%%%%%%%%%%%%%

\section{Wide and narrow ribbons of nematic elastomers}
\label{sec:main}

%%%%%%%%%%%%%%%%%%%%%%%%%%%%%%%%%%%%%%%%%%
%%%%%%%%%%%%%%%%%%%%%%%%%%%%%%%%%%%%%%%%%%

In this section, we consider a thin sheet 
of nematic elastomer with \emph{twist} 
distribution of the director along the thickness.
Starting from a $3D$ model, described in 
Subsection \ref{3D}, 
we present in Subsection \ref{2D}
the rigorous derivation of a corresponding
$2D$ plate model. In Subsection \ref{1D}
we then derive, again via $\Gamma$-convergence arguments,
a $1D$ rod model from the $2D$ model.

Before proceeding, let us establish some general
notation which will be used throughout.
For the standard basis of $\R^3$ and $\R^2$,
we use the notation
$\{\mathsf e_1,\mathsf e_2,\mathsf e_3\}$
and
$\{\mathsf f_1,\mathsf f_2\}$,
respectively.
$\SO(3)$ and ${\rm Sym}(3)$ are the sets of $3{\times3}$ 
rotations and symmetric matrixes, respectively,
and ${\rm I}\in\SO(3)$ is the identity matrix.
We use the symbol $\tr^2A$ for the square
of the trace of a matrix $A$,
and denote by $\Sph^2$ the unit sphere of $\R^3$.
Finally, the symbol $C$ will be used
to denote a positive constant depending on given data, 
and possibly varying from line to line.

%%%%%%%%%%%%%%%%%%%%%%%%%%%%%%%

\subsection{The three-dimensional model}
\label{3D}

%%%%%%%%%%%%%%%%%%%%%%%%%%%%%%%

\begin{figure}[htbp]
\begin{center}
\includegraphics[width=7cm]
{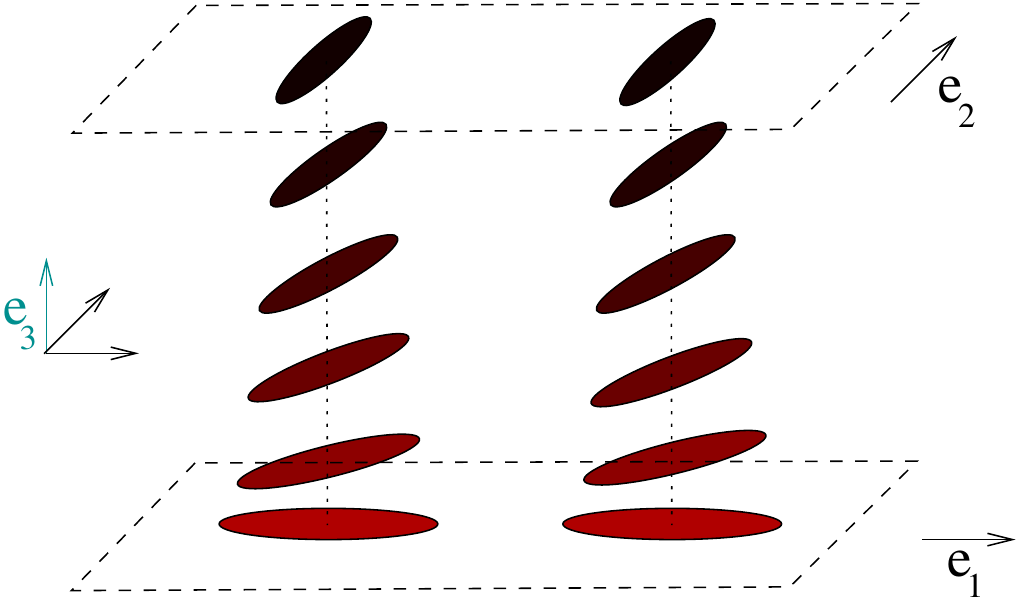}
\end{center}
\caption{
Sketch of the twist director field.
} 
\label{fig:twist}
\end{figure}

We consider a thin sheet of nematic elastomer occupying the
reference configuration 
\begin{equation*}
%\label{def:ref_config}
\Om_h^\ep
\,:=\,
\Om_h^\ep(\theta)
\,:=\,
\bigg\{
z_1\mathsf e_1^{\theta}+z_2\mathsf e_2^{\theta}
+z_3\mathsf e_3\,:\,
z_1\in\big(-\ell/2,\ell/2\big),\ z_2\in\big(-\ep/2,\ep/2\big),
\ z_3\in\big(-h/2,h/2\big)
\bigg\},
\end{equation*}
where
\[
\mathsf e_1^{\theta} 
\,:=\, 
\left(
\begin{array}{c}
\cos\theta\\
\sin\theta\\
0
\end{array}
\right)
\qquad\quad\mbox{and}\qquad\quad
\mathsf e_2^{\theta} 
\,:=\, 
\left(
\begin{array}{c}
-\sin\theta\\
\cos\theta\\
0
\end{array}
\right),
\]
and $\ep>0$ and $h>0$ are small dimensionless parameters
such that
\[
\ell\gg\ep\gg h.
\]
The directions $\mathsf e_1$
and $\mathsf e_2$ play a special role,
which is related to the orientation of the
nematic director on the top and bottom faces of the thin sheet,
see \eqref{twist_geom} below and Figure \ref{fig:twist}.
For future reference, we also set
$\Om_h^\ep
:=
\omega^\ep\times(-h/2,h/2)$,
with
\begin{equation}
\label{omega_ep}
\omega^\ep
\,=\,
\omega^\ep(\theta)
\,:=\,
\bigg\{
z_1\mathsf f_1^{\theta}+z_2\mathsf f_2^{\theta}
\,:\,
z_1\in\big(-\ell/2,\ell/2\big),\ z_2\in\big(-\ep/2,\ep/2\big)
\bigg\},
\end{equation}
and
\begin{equation}
\label{cose}
\mathsf f_i^{\theta}
\,:=\,
R_\theta\mathsf f_i,\quad\ i=1,2,
\qquad\quad
R_\theta 
\,:=\, 
\left( 
\begin{array}{cc} 
\cos\theta & -\sin\theta  \\ 
\sin\theta & \cos\theta  \\
\end{array}
\right),
\qquad\quad
0\leq\theta<\pi.
\end{equation}
Our nematic elastomer thin sheet will be modelled by   
a stored energy density
which is heterogeneous along the thickness.
More precisely, 
in the model we are going to consider,
the nematic director 
varies along the thickness. This
will induce a $z_3$-dependence of
the spontaneous strain distributions
and, in turn,
$z_3$-dependent stored energy densities.
The system under analysis is a \emph{twist} 
nematic elastomer thin sheet. 
The distribution of the nematic director in 
the twist geometry is defined as
\begin{equation}
\label{twist_geom}
\hat n_h(z_3)
\,:=\,
\left(
\begin{array}{c}
\cos\big(\frac{\pi}4+\frac\pi2\frac{z_3}h\big)\\
\sin\big(\frac{\pi}4+\frac\pi2\frac{z_3}h\big)\\
0
\end{array}
\right),
\qquad\qquad
z_3\in\Big[-\frac h2,\frac h2\,\Big],
\end{equation}
see Figure \ref{fig:twist}.
Note that this distribution is constant
on each horizontal plane and in particular $n_h$ is
(constantly) equal to $\mathsf e_1$ and to to $\mathsf e_2$
on the bottom and on the top face of the
sheet, respectively.   
Also, recall that such distribution is the 
($\ep$-independent) solution to the minimum problem
\begin{equation*}
\min_{\begin{array}{r}n(-h/2)=\mathsf e_1\\
n(h/2)=\mathsf e_2
\end{array}}
\int\limits_{\Omega_h^\ep}|\na n|^2{\rm d}z.
\end{equation*}
This is the orientation arising in the fabrication procedure
of the material, in which the director is oriented in the
liquid phase, and then ``frozen'' by the 
crosslinking process transforming the liquid into an elastomer.

Now, if $n\in\R^3$ is a unit vector representing 
the local order of the nematic director, 
the (local) response of the nematic elastomer is
encoded by the positive definite symmetric matrix
\begin{equation}
\label{nem_tens}
L(n)=a^{\frac23}n\otimes n+a^{-\frac13}({\rm I}-n\otimes n),
\end{equation} 
where $a>1$ is a material parameter.
More precisely, under the effect of 
lowering the temperature below a threshold
(isotropic-to-nematic transition temperature),
in a region where the nematic director
is $n$, the material spontaneously deforms
via a deformation $y$ such that
$\na y^T\na y~=~L(n)$.
As a consequence, 
the through-the-thickness variation of the nematic director 
translates into the spontaneous strain field 
\begin{equation}
\label{spont_phys}
\widehat C_h(z_3)
\,:=\,
L(\hat n_h(z_3))
\,=\,
a_h^{2/3}\hat n_h(z_3)\otimes\hat n_h(z_3)
+a_h^{-1/3}\big({\rm I}-\hat n_h(z_3)\otimes\hat n_h(z_3)\big).
\end{equation}
Notice that, in this expression, we allow the material 
parameter $a$ to be $h$-dependent.
More precisely, from now on we will assume that  
\begin{equation}
\label{a_h}
a_h=1+\frac{\alpha_0}{h_0}h,
\end{equation}
where $\alpha_0$ is a positive dimensionless parameter, while $h_0$
and $h$ have the physical dimension of length.
This assumption is easily understandable
if one thinks that,
arguing as in the proof of
the Gauss Lemma
(see, e.g., \cite[Chapter 5]{Petersen}), 
a given metric $g$ of $\R^3$
can always be expressed around
a smooth hypersurface $\Sigma$ as
\[
g
\,=\,
dr\otimes dr+g_{ij}d\theta^i\otimes d\theta^j,
\]
using coordinates $\{r,\theta^1,\theta^2\}$, 
where $r$ is the signed $g$-distance from $\Sigma$.
In the coordinates $\{\theta^1,\theta^2\}$,
one can then express the 
second fundamental form ${\rm A}$ on $\Sigma$ as
\[
{\rm A}_{ij}\,d\theta^i\otimes d\theta^j
\,=\,
\frac12\frac{\pa g_{ij}}{\pa r}\,d\theta^i\otimes d\theta^j.
\] 
In other words, via the standard identification of
metric and strain, curvature is related to the ratio
between the magnitude of the strain difference along the thickness and the thickness itself. 
Hence, the linear scaling in $h$ in \eqref{a_h} is needed
in order to obtain a finite curvature in the limit $h\to0$.

To model our system in the framework of finite elasticity,
we consider the energy density 
$\widehat W_h$ defined on points of 
$(-h/2,h/2){\times}\RRR$ as 
\begin{equation}
\label{prot_energy}
\widehat W_h(z_3,F):=
\left\{
\begin{array}{ll}
\displaystyle
\frac{\mu}2\Big[(F^TF)\cdot\widehat C_h^{-1}(z_3)
-3-2\log(\det F)\Big]
+W_{vol}(\det F) & \quad\mbox{if }\det F>0,\\
+\infty & \quad\mbox{if }\det F\leq0,
\end{array}
\right.
\end{equation}
which also accounts for compressibility effects.
Here, $\mu>0$ is a material constant 
(shear modulus) and the function
$W_{vol}:(0,\infty)\to[0,\infty)$ is ${\rm C}^2$ around $1$
and fulfills the conditions:
\begin{equation*}
W_{vol}(t)=0 \iff t=1,
\qquad
W_{vol}(t)\longrightarrow\infty\ \,{\rm as}\ \,t\to 0^+,
%\ \
%W_{vol}(t) \geq C\,t^2\ \ \mbox{for }t>>0,
\qquad
W_{vol}''(1) >0.
\end{equation*}
A typical example of $W_{vol}$ is
$t\mapsto c\,(t^2-1-2\log t)$, 
where $t$ is nondimensional and $c$, a positive constant,
is energy per unit volume, 
so that $W_{vol}''(1)$ has the same dimension
of $\mu$.
Expression \eqref{prot_energy} is
a natural generalization, see \cite{Ag_DeS_1}, 
of the classical trace formula for nematic elastomers 
derived by Bladon, Terentjev and Warner \cite{Bla_Ter_War},
in the spirit of Flory's work 
%pagg 468-469 Flory 
on polymer elasticity \cite{Flo}.
The presence of the purely volumetric term $W_{vol}(\det F)$
guarantees that the Taylor expansion at order two of the density 
results in isotropic elasticity with two independent 
elastic constants (shear modulus and bulk modulus).

Observe that, for every $F$ with $\det F>0$,
\begin{equation}
\label{W_h_W_0}
\widehat W_h(z_3,F)
\,=\,
W_0\big(F\widehat U_h(z_3)^{-1}\big),
\qquad
\widehat U_h(z_3)
\,:=\,
\widehat C_h(z_3)^{\frac12},
\end{equation}
where
\begin{equation}
\label{nostra_W_0}
W_0(F)
\,:=\,
\frac{\mu}2
\Big[
|F|^2-3-2\log(\det F\,)
\Big]
+
W_{vol}(\det F\,).
\end{equation}
Also, note that 
\[
W_0(F)
\,\geq\,
\frac{\mu}2
\Big[
|F|^2-3-2\log(\det F\,)
\Big]
\,\geq\,
\frac32\mu
\big[
(\det F)^{\frac23}-1-\log(\det F)^{\frac23}
\big]
\,\geq\,
0,
\]
where in the second inequality we have used that
$|F|^2\geq 3(\det F)^{\frac23}$, due to the inequality
between arithmetic and geometric mean,
and the third inequality comes from direct computations.
Also, the second and the third inequalities become 
equalities iff $F^TF=\alpha{\rm I}$, for some $\alpha\geq0$,
and iff $\det F=1$, respectively.
All in all, we have that $W_0$ is a nonnegative function
minimised precisely at ${\rm SO}(3)$.
In turn, from \eqref{W_h_W_0}, we have that
$\widehat W_h(z_3,\cdot)$ is a nonnegative function
minimised precisely at ${\rm SO}(3)\widehat C(z_3)^{\frac12}$.
Let us also point out that from the expression
of the second differential ${\rm D}^2W_0({\rm I})$  
of $W_0$ at ${\rm I}$ applied twice to some 
$M\in\R^{3\times3}$, namely
\begin{equation}
\label{diff_secondo}
Q_3(M)
\,:=\,
{\rm D}^2W_0({\rm I})[M,M]
\,=\,
2\mu\,|\sym M|^2+W_{vol}''(1)\,\tr^2M,
\end{equation}
and from the fact that $W_0(F)\geq C|F|^2$ for $|F|\gg0$,
one can easily deduce that
\begin{equation}
\label{crescita}
W_0(F)
\,\geq\,
C{\rm dist}^2(F,{\rm SO}(3)),
\qquad\quad
\mbox{for all}\ F\in\R^3.
\end{equation}
This is one of the sufficient conditions
which allow us to apply the results of \cite{Schmidt2007}, 
in the next subsection.
Note that the quadratic growth of $W_0$
from ${\rm SO}(3)$ implies that
\begin{align*}
\widehat W_h(z_3,F)
\,\geq\,
C{\rm dist}^2\big(F\widehat U_h(z_3)^{-1},{\rm SO}(3)\big)
\,\geq\,
C\lambda_{min}^2(z_3)
{\rm dist}^2\big(F,{\rm SO}(3)\widehat U_h(z_3)\big),
\qquad\quad
\mbox{for all}\ F\in\R^3,
\end{align*}
where $\lambda_{min}(z_3)>0$ is the minimum eigenvalue
of $\widehat U_h(z_3)$.
Also, observe from \eqref{diff_secondo} that
while $\mu$ is the elastic shear modulus,
$W_{vol}''(1)+2\mu/3$ has the physical meaning 
of bulk modulus.

We denote by $\widehat{\mathscr F}_h^\ep$ the free-energy
functional corresponding to the energy density $\widehat W_h$, i.e. 
\begin{equation}
\label{phys_quan}
\widehat{\mathscr F}_h^\ep(v)
\,:=\,
\int\limits_{\Om_h^\ep}\widehat W_h(z_3,\na v(z))\,{\rm d}z,
\end{equation}
where $v:\Om_h^\ep\to\R^3$ is a deformation
and $z=(z',z_3)$ is a point of $\Om_h^\ep$, 
whose components are given w.r.t.~ the canonical
basis $\{\mathsf e_1,\mathsf e_2,\mathsf e_3\}$.

%%%%%%%%%%%%%%%%%%%%%%%%%%%%%%%

\subsection{The plate model}
\label{2D}

%%%%%%%%%%%%%%%%%%%%%%%%%%%%%%%

In order to  derive rigorously, in the limit as $h\downarrow0$, 
a 2D model from the 3D setting, 
we first make a standard change of variables which
allows us to rewrite the energies in a fixed, $h$-independent 
rescaled reference configuration.
We denote by $x=(x_1,x_2,x_3)=(x',x_3)$ 
an arbitrary point in the rescaled reference configuration 
$\Om^\ep:=\omega^\ep\times(-1/2,1/2)$.
For every $h>0$ small, we define the rescaled energy density
$W_h:(-1/2,1/2)\times\R^{3\times3}\longrightarrow[0,+\infty]$
%and the rescaled applied loads
%$f_h^\ep:\Om^\ep\longrightarrow\R^3$ 
as
\begin{equation}
\label{eq:resc_ener}
W_h(x_3,F)
\,:=\,
\widehat W_h(hx_3,F).
%\qquad\qquad\quad
%f_h^\ep(x)
%\,:=\,
%\hat f_h^\ep(x',hx_3).
\end{equation}
Note that $W_h$ fulfills
\begin{equation*}
W_h(x_3,F)=0\qquad\quad\mbox{iff}
\qquad\quad
F\in{\rm SO}(3)\overline U_h(x_3),
\qquad\qquad
\overline U_h(x_3):=\widehat U_h(hx_3).
\end{equation*}
We also set $\overline C_h(x_3):=\overline U_h^2(x_3)$.
Finally, defining 
\begin{equation}
\label{resc_grad}
\na_hy:=
\left(
\begin{array}{ccccc}
\!\!\partial_{x_1}y\!\! & \!\!\bigg|\!\! 
& \!\!\partial_{x_2}y\!\! & \!\!\bigg|\!\! 
& \!\!\displaystyle\frac{\partial_{x_3}y}h\!\!
\end{array}
\right)
=:
\left(
\begin{array}{ccccc}
\!\!\na'y\!\!
& \!\!\bigg|\!\! 
& \!\!\displaystyle\frac{\partial_{x_3}y}h\!\!
\end{array}
\right),
\qquad\qquad
\mbox{for every}\quad y:\Om^\ep\longrightarrow\R^3,
\end{equation}
the correspondence between the original quantities and the
rescaled ones is through the formulas
\begin{equation}
\label{correspondence}
%\hat{\mathscr E}^h(v)
%\,=\,
%h\,\mathscr E^h(y),
%\qquad\quad
\widehat{\mathscr F}_h^\ep(v)
\,=\,
h\,\mathscr F_h^\ep(y),
\qquad\qquad
v(z):=y(z',z_3/h)
\quad\mbox{for a.e.}\ z\in\Omega_h^\ep.
\end{equation}
Here, 
% the rescaled stored elastic energy functional $\mathscr E^h$ 
the rescaled free-energy functional $\mathscr F_h^\ep$
is defined, 
on a deformation $y:\Om^\ep\to\R$, as
\begin{equation*}
%\label{mathscr_E_h}
\mathscr F_h^\ep(y)
\,:=\,
\int\limits_{\Om^\ep}W_h(x_3,\na_h y(x))\,\dd x.
%-
%\int\limits_{\Om^\ep}f_h^\ep\cdot y\,\dd x.
\end{equation*}

To deal with the (rescaled)
spontaneous strains $\overline C_h(x_3)$,
we first note that the rescaled director field
has the following expression:
\begin{equation}
\label{tensor_N}
n(x_3)
:=
\hat n_h(hx_3)
=
\left(
\begin{array}{c}
\cos\big(\frac{\pi}4+\frac\pi2x_3\big)\\
\sin\big(\frac{\pi}4+\frac\pi2x_3\big)\\
0
\end{array}
\right)\!,
\end{equation}
which is independent of $h$.
In turn, we obtain for $\overline C_h(x_3)$
the equivalent expression
\begin{align}
\label{eq:resc_SS}
\overline C_h(x_3)
&
=a_h^{2/3}n(x_3)\otimes n(x_3)
+a_h^{-1/3}\big({\rm I}-n(x_3)\otimes n(x_3)\big)\nonumber\\
&
=\big(a_h^{2/3}-a_h^{-1/3}\big)
\left(\frac{\rm I}{a_h-1}+n(x_3)\otimes n(x_3)\right)\nonumber\\
&
={\rm I}-2\,h\,B(x_3)+R^h(x_3),
\qquad\quad
B(x_3)
\,:=\,
\frac12\frac{\alpha_0}{h_0}
\left(\frac{\rm I}3-n(x_3)\otimes n(x_3)\right),
\end{align}
where 
$\|R^h\|_{\infty}=o(h)$
and $\|\cdot\|_{\infty}$ is the norm in the space
${\rm L}^{\infty}\big((-1/2,1/2),\R^{3\times3}\big)$.
Note that in the third equality we have
plugged in expression \eqref{a_h} for $a_h$
and used the expansion
\begin{equation*}
a_h^{2/3}-a_h^{-1/3}
=
\frac{\alpha_0}{h_0}h
-\frac13\left(\frac{\alpha_0}{h_0}\right)^2h^2
+o(h^2).
\end{equation*}
Moreover, plugging into
\eqref{eq:resc_SS} the expression for $n$ given by
\eqref{tensor_N}, we have that
\begin{equation}
\label{B_expl}
B(x_3)
\,=\,
\frac12\frac{\alpha_0}{h_0}
\left(
\begin{array}{ccc}
\frac13-\cos^2\big(\frac{\pi}4+\frac\pi2x_3\big) 
& \frac12\sin\big(\frac{\pi}2+\pi x_3\big) & 0 \\
\frac12\sin\big(\frac{\pi}2+\pi x_3\big) 
&  \frac13- \sin^2\big(\frac{\pi}4+\frac\pi2x_3\big) & 0 \\
0 & 0 & \frac13\\
\end{array}
\right).
\end{equation}
The description of the 3D system provided by
\eqref{W_h_W_0}, \eqref{eq:resc_ener}, 
and \eqref{eq:resc_SS} allows us to
directly use the results of \cite{Schmidt2007}
for the derivation of a corresponding limiting model,
in the regime of finite bending energies.
More precisely, in \cite{Ag_De_bend} 
we can state a compactness result for sequences 
$\{y_h\}$ such that $\mathscr F_h^\ep(y_h)\leq Ch^2$,
and a $\Gamma$-convergence result for 
$\mathscr F_h^\ep/h^2$, in the limit $h\downarrow0$.
In particular, we can compute the limiting 2D model.
In order to do it, we first compute,
for every $G\in\R^{2\times 2}$, 
the relaxed energy density 
\begin{equation}
\label{Q2}
Q_2(G)
\,:=\,
\min_{\stackrel{b\in\R^2}{a\in\R}}
Q_3
\left(
\left[
\begin{tabular}{c|c}
$G$ & $b$\\
\hline
$0$ & $a$
\end{tabular}
\right]
\right)
\,=\,
2\mu\Big(\,|\sym G|^2+\gamma\,\tr^2G\Big),
\end{equation}
where $Q_3$ is defined in \eqref{diff_secondo} and
\begin{equation}
\label{def:gamma}
\gamma:=
\frac{W_{vol}''(1)}{2\mu+W_{vol}''(1)}.
\end{equation}
The adjective ``relaxed'' used for the 
two-dimensional energy density 
$Q_2$ is due to the fact that 
it arises from an optimisation procedure.
Then, we compute the doubly relaxed energy density
\begin{align}
\label{barQ2}
\overline Q_2(G)
&\,:=\,
\min_{D\in\R^{2\times2}}\int\limits_{-1/2}^{1/2}
Q_2\big(D+t\,G+\check B(t)\big){\rm d}t\\
&\,=\,
\alpha_T\,Q_2(G-\bar A_T)+\beta_T,\nonumber
\end{align}
where the $2{\times}2$ symmetric matrix $\check B$ 
is obtained from $B$ 
(cf. \eqref{eq:resc_SS} and \eqref{B_expl}) 
by omitting the last row and the last column,
namely
\begin{equation}
\label{cek_B}
\check B
\,:=\,
\sum_{ij=1}^2B_{ij}\,\mathsf f_i\otimes\mathsf f_j.
\end{equation}
Moreover, 
the constants $\alpha_T$
(a geometric parameter, reminiscent of the 
``moment of inertia'' of a cross section of unit width), 
$\bar A_T$ (the \emph{target curvature} tensor),
and $\beta_T$ (a positive constant, reminder of the 
presence of residual stresses, cf. Remark \ref{frustration})
are given by the formulas:
\begin{equation}
\label{const_T}
\alpha_T\,=\,\frac1{12},
\quad\qquad
\bar A_T\,=\,\frac{6}{\pi^2}\frac{\alpha_0}{h_0}
\,{\rm diag(-1,1)},
\quad\qquad
\beta_T\,=\,\mu\,\Big(\frac{\pi^4-4\pi^2-48}{4\pi^4}\Big)
\frac{\alpha_0^2}{h_0^2}.
\end{equation}
The 2D free-energy limit functional
turns out to have the following expression:
\begin{equation}
\label{twist_2D}
\mathscr F_{lim}^\ep(y)
\,:=\,
\frac12
\int\limits_{\omega^\ep}
\left[
\alpha_T\,Q_2(A_y(x')-\bar A_T)+\beta_T
\right]
{\rm d}x',
\end{equation}
if $y\in{\rm W}^{2,2}_{\rm iso}(\omega^\ep,\R^3)$,
while 
$\mathscr F_{lim}^\ep(y)=+\infty$ 
elsewhere in
${\rm W}^{1,2}(\Omega^\ep,\R^3)$.
Here, 
the symbol $A_y(x')$ stands for the second fundamental form
associated with $y(\omega_\ep)$ and defined at points $x'\in\omega_\ep$,
namely $A_y(x'):=(\na'y(x'))^T\na'\nu(x')$, with
$\nu:=\pa_1y\wedge\pa_2y$.
Moreover,
the class ${\rm W}^{2,2}_{\rm iso}(\omega^\ep,\R^3)$
is that of the isometric immersions of $\omega^\ep$ into the three-dimensional
Euclidean space, namely
\begin{equation}
\label{Aiso}
{\rm W}^{2,2}_{\rm iso}(\omega^\ep,\R^3)
\,:=\,
\Big\{
y\in{\rm W}^{2,2}(\omega^\ep,\R^3)\,:\,(\na'y)^T\na'y={\rm I}_2
\ \mbox\ {a.e.\ in\ }\,\omega^\ep
\Big\},
\end{equation}
where ${\rm I}_2$ is the identity matrix of $\R^{2\times2}$.
%Also, the ``reduced'' loading $f_{lim}^\ep$ is defined
%in terms of the limiting loading $f^\ep$ as 
%\begin{equation}
%\label{f_lim}
%f_{lim}^\ep(x')
%\,:=\,
%\int\limits_{-1/2}^{1/2}
%\!\!f^\ep(x',x_3)\,\dd x_3,
%\qquad\quad\mbox{for a.e.}\quad x'\in\omega^\ep.
%\end{equation}

The following theorem is a straightforward corollary of the aforementioned
compactness and $\Gamma$-convergence results, for which we 
refer the reader to \cite{Ag_De_bend}.

\begin{theorem}
\label{main_thm_T}
Setting
$$
m_h:=\inf_{y\in{\rm W}^{1,2}}\mathscr F_h^\ep,
$$
suppose that $\{y_h\}$ is a low-energy sequence, i.e.
it fulfills
\begin{equation*}
%\label{almost_min_T}
\lim_{h\to0}\frac{\mathscr F_h^\ep(y_h)}{h^2}
\,=\,
\lim_{h\to0}\frac{m_h}{h^2}.
\end{equation*}
Then, up to a subsequence, $y_h\longrightarrow y$ in 
${\rm W}^{1,2}(\Om^\ep,\R^3)$, where
$y\in{\rm W}^{2,2}_{\rm iso}(\omega^\ep,\R^3)$ 
is a solution to the minimum problem
\begin{equation*}
m_{lim}:=\min_{{\rm W}^{2,2}_{\rm iso}(\omega^\ep,\R^3)}
\mathscr F_{lim}^\ep.
\end{equation*}
Moreover, $(m_h/h^2)\to m_{lim}$.
\end{theorem}

Let us now fix a low-energy sequence $\{y_h\}$ 
converging to a minimiser 
$y\in{\rm W}^{2,2}_{\rm iso}(\omega^\ep,\R^3)$
of the 2D model \eqref{twist_2D},
and rephrase the theorem in terms of the physical
quantities \eqref{phys_quan}. 
Defining the deformations $v_h(z',z_3)=y_h(z',z_3/h)$
in the physical reference configuration $\Omega_h$,
we have
$\lim_{h\to0}\hat{\mathscr F}_h^\ep(v_h)/h^3
=
m_{lim}$,
in view of \eqref{correspondence}.
Equivalently, for a given small thickness $h_0$,
the approximate identity
\begin{equation*}
\hat{\mathscr F}_T^{h_0}(v_{h_0})
\,\cong\,
\frac{\mu\,h_0^3}{12}\int\limits_{\omega^\ep}
\left\{\bigg|{\rm A}_y(x')-
\frac{6}{\pi^2}\frac{\alpha_0}{h_0}
\,{\rm diag(-1,1)}
\bigg|^2
+\gamma\,
{\rm H}_y^2(x')\right\}\dd x'
\,+\, \mu\,h_0\,\alpha_0^2
\Big(\frac{\pi^4-4\pi^2-48}{8\pi^4}\Big)|\omega^\ep|
\end{equation*}
holds true,
modulo terms of order higher than $3$ in $h_0$.
Here, the symbol ${\rm H}_y$ 
denotes the mean curvature of 
$y(\omega^\ep)$,
namely ${\rm H}_y=\tr A_y$.
In reading the formula, 
recall that the second fundamental form has 
physical dimension of inverse length.

\begin{remark}
[Kinematic incompatibility in 3D and 
geometric obstructions in 2D]
\label{frustration}
We remark that the spontaneous strain field
$\Omega_h^\ep\ni z\mapsto \widehat C(z)=\widehat C(z_3)$ 
defined in \eqref{spont_phys}
is not \emph{kinematically compatible}, meaning
that there are no smooth deformations
$v:\Omega_h^\ep\to\R^3$ such that, for every $z\in\Omega_h^\ep$, 
$\det\na v(z)>0$ and
\begin{equation*}
%\label{kin_com}
\na v^T\na v(z)=\widehat C(z_3).
\end{equation*}
To prove this fact, one can interpret $\widehat C$ as
a given metric on $\Omega_h^\ep$ and show equivalently
(see \cite{Ag_De_bend} and \cite{Ciarlet_book})
that the fourth-order Riemann curvature tensor
associated with the metric $\widehat C$ is 
not identically zero in $\Omega_h^\ep$.
This kinematic incompatibility and the related 
impossibility of minimising (to the value zero) the 3D 
energy functional \eqref{phys_quan},
results into the presence of the positive constant
$\beta_T$ in the 2D limiting model \eqref{const_T}--\eqref{twist_2D},
which reveals the presence of residual stresses.
Note at the same time that the quadratic term 
$Q_2(\cdot-\bar A_T)$ in \eqref{twist_2D} per se 
prevents the attainment of zero as a minimum value of the energy.
Indeed, there are no deformations $y$ in the class
${\rm W}^{2,2}_{\rm iso}(\omega^\ep,\R^3)$ such that
$y(\omega^\ep)$ has non-zero Gaussian curvature
(while $\det \bar A_T<0$).
More precisely, one can prove
(see \cite[Lemma 3.8]{Ag_De_bend})
that
\begin{equation*}
\min_{{\rm W}^{2,2}_{\rm iso}(\omega^\ep,\R^3)}
\mathscr F^{lim}
\,=\,
\mathscr F^{lim}(y_T)
\,=\,
\frac{3\,\mu}{\pi^4}\frac{\alpha_0^2}{h_0^2}
\left(\frac{1+2\gamma}{1+\gamma}\right)|\omega^\ep|
+\frac{\beta_T}2|\omega^\ep|,
\end{equation*}
where $y_T\in{\rm W}^{2,2}_{\rm iso}(\omega^\ep,\R^3)$
is such that
\begin{equation}
\label{cond_min}
\mbox{either}\qquad
{\rm A}_{y_T}
\equiv
\frac6{\pi^2(1+\gamma)}\frac{\alpha_0}{h_0}
{\rm diag}\left(-1,0\right)
\qquad\mbox{or}\qquad
{\rm A}_{y_T}
\equiv
\frac6{\pi^2(1+\gamma)}\frac{\alpha_0}{h_0}
{\rm diag}\left(0,1\right).
\end{equation}
Recall that $\omega^\ep=\omega^\ep(\theta)$
and therefore the minimisers $y_T$ are $\theta$-dependent
in the sense that the long axis of the sheet
($\mathsf f_1^\theta$ in \eqref{omega_ep})
makes angles $\theta$ and $\pi/2-\theta$ 
with the eigenvectors 
$\mathsf f_1$ and $\mathsf f_2$ corresponding to 
the nonzero eigenvalue of the realised second fundamental
forms \eqref{cond_min},
respectively.

\end{remark}

%%%%%%%%%%%%%%%%%%%%%%%%%%%%%%%

\subsection{The rod model}
\label{1D}

%%%%%%%%%%%%%%%%%%%%%%%%%%%%%%%

In this subsection, following \cite{Ago_DeS_Kou}, 
we derive a 1D reduced model, 
in the limit of vanishing width $\ep\downarrow0$, 
starting from the 2D setting 
\eqref{const_T}--\eqref{twist_2D}.
More precisely, the starting
point of the following 1D derivation 
is the 2D bending
energy obtained by multiplying expression \eqref{twist_2D} 
by the factor $\ep^{-1}$, namely
\begin{align*}
\widehat{\mathscr F}^{\ep}(\hat y)
\,:=\,
\frac1{\ep}
\mathscr F_{lim}^\ep(\hat y)
&\,=\,
\frac1{2\,\ep}
\int\limits_{\omega^\ep}
\Big\{
\alpha_T\,Q_2(A_{\hat y}(x')-\bar A_T)+\beta_T
\Big\}
{\rm d}x'\\
&\,=\,
\frac1{\ep}
\int\limits_{\omega^\ep}
\left\{
\frac{\mu}{12}
\left|A_{\hat y}(x')-\bar A_T\right|^2+
\frac{\gamma\mu}{12}\,
{\rm H}_{\hat y}^2(x')
+\frac{\beta_T}2
\right\}
{\rm d}x',
\end{align*}
where $\hat y$ is a deformation in the class
${\rm W}^{2,2}_{\rm iso}(\omega^\ep,\R^3)$,
and $\bar A_T$ and $\beta_T$ are given in \eqref{const_T}.
Now, observe that the functional
$\widehat{\mathscr F}^{\ep}$ is $\theta$-dependent,
$\theta\in[0,\pi)$, 
since our thin and narrow nematic elastomer sheet
has been ``cut out'' from the plane 
$(\mathsf e_1,\mathsf e_2)$
at an angle $\theta$ with the horizontal axis
(see definition \eqref{omega_ep} of $\omega^\ep$).
To simplify the notation, but also to emphasize the
$\theta$-dependence of our model, we perform a change
of variable from the $\theta$-rotated strip 
$\omega^\ep(\theta)$ to the horizontal strip 
$\omega^\ep(0)$ and introduce,
for a deformation $\hat u$ defined on $\omega^\ep(0)$,
the energy
\begin{equation*}
%\label{rescaled_energy}
\widehat{\mathscr E}^{\ep}(\hat u)
\,:=\,
\frac1{\ep}
\int\limits_{\omega^\ep(0)}
\left\{
\frac{\mu}{12}
\left|A_{\hat y}(x')-\bar A_T^\theta\right|^2+
\frac{\gamma\,\mu}{12}\,
{\rm H}_{\hat y}^2(x')
+\frac{\beta_T}2
\right\}
{\rm d}x',
\end{equation*}
where
\begin{equation*}
%\label{A_theta}
\bar A^{\theta}_T
\,:=\,R_\theta^{\rm T}\bar A_TR_\theta,
\end{equation*}
and $R_\theta$ is defined as in \eqref{cose}.
Note that for every $\theta\in[0,\pi)$,
we have
\begin{equation}
\label{A^T}
\bar A_T^\theta 
\,=\,
\kk\left( \begin{array}{rr} -a_\theta & b_\theta \\ 
b_\theta & a_\theta \end{array}\right),
\qquad\ 
a_\theta:=\cos2\theta,
\ \quad
b_\theta:=\sin2\theta,
\ \quad
\kk:=\frac 6{\pi^2}\frac{\alpha_0}{h_0}.
\end{equation}
Some computations show that indeed 
$\widehat{\mathscr E}^{\ep}(\hat u)
=
\widehat{\mathscr F}^{\ep}(\hat y),
$
whenever 
$\hat u(x')=\hat y(R_\theta x')$,
$x'\in\omega^\ep(0)$.
Keeping this identification in mind,
we now proceed similarly to the 3D-to-2D derivation
and express the above energy over the fixed 2D domain
\[
\omega
\,:=\,
I\times\Big(\!-\frac12,\frac12\,\Big),
\qquad\qquad
I := \Big(\!-\frac\ell2,\frac\ell2\,\Big),
\]
In order to do this, 
we operate another change variables and define for a 
deformation $\hat u:\omega^\ep(0)\to\R^3$ 
its rescaled version  
$u:\omega\to\R^3$, given by
\[
u(x_1,x_2) 
\,=\, 
\hat u(x_1,\ep\,x_2),
\qquad\qquad
(x_1,x_2)\in\omega
\]
(note that we use the same notation $x'=(x_1,x_2)$
for points belonging to the different sets
$\omega^\ep(\theta)$, $\omega^\ep(0)$, and $\omega$,
when there is no ambiguity).
Moreover, introducing the scaled gradient
\[
\nabla_\ep \cdot = (\partial_1 \cdot |\ep^{-1}\partial_2\cdot)
\]
we obtain that $\nabla_\ep u(s,t) = \nabla\hat u(s,\ep t)$ 
and $u$ belongs to the space of scaled isometries of 
$\omega$ defined as
\[
{\rm W}^{2,2}_{{\rm iso},\ep}(\omega,\R^3)
\,:=\,
\Big\{u\in{\rm W}^{2,2}(\omega,\R^3)\,:\, 
|\partial_1 u| = |\ep^{-1}\partial_2 u|=1,\,\,
\partial_1 u\cdot\partial_2 u = 0\mbox{ a.e. in }
\omega\Big\}.
\]
Similarly, we may define the scaled unit normal to 
$u(\omega)$ by
\[
n_{u,\ep} = \partial_1 u\wedge\ep^{-1}\partial_2 u
\]
and the scaled second fundamental form associated 
with $u(\omega)$ by
\[
A_{u,\ep} =  \left( \begin{array}{rrr} n_{u,\ep}\cdot \partial_1\partial_1 u & \, & \ep^{-1}n_{u,\ep}\cdot \partial_1\partial_2 u \\  \ep^{-1}n_{u,\ep}\cdot \partial_1\partial_2 u &\, &  \ep^{-2}n_{,\ep}\cdot \partial_2\partial_2 u \end{array}\right).
\]
With this definition, 
$A_{y,\ep}(x_1,x_2) = A_v(x_1,\ep x_2)$, 
and 
$\widehat{\mathscr E}^\ep(\hat u) 
= 
\mathscr E^\ep(u)$, 
where 
\begin{equation*}
%\label{our_fun}
\mathscr E^\ep(u) 
\,:=\, 
\int\limits_{\omega}
\left\{
\frac{\mu}{12}
\left|A_{u,\ep}(x')-\bar A_T^\theta\right|^2+
\frac{\gamma\mu}{12}\,
\tr^2A_{u,\ep}(x')
+\frac{\beta_T}2
\right\}
{\rm d}x',
\qquad\quad
\mbox{for every }u\in
{\rm W}^{2,2}_{{\rm iso},\ep}(\omega,\R^3).
\end{equation*}

An adaptation of the theoretical results of 
\cite{Freddi2015} allows to prove in \cite{Ago_DeS_Kou}
a compactness result for sequences 
$\{u^\ep\}$ such that $\mathscr E^\ep(u^\ep)$ 
is uniformly bounded, and a $\Gamma$-convergence result
for the functionals $\mathscr E^\ep$.
We gather in Theorem \ref{cor:twist} below 
the most important consequences
of these results. 
The limiting 1D free-energy functional 
$\mathscr E^{lim}$
turns out to have the
following expression:
\begin{equation}
\label{funct_lim_T}
\mathscr E^{lim}(d_1,d_2,d_3)
\,:=\,
\int_I 
\overline Q_T^{\,\theta}
(d_1^\prime\cdot d_3,d_2^\prime\cdot d_3)\,{\rm d}x_1,
\end{equation}
on every triplet $(d_1,d_2,d_3)$, representing
a orthonormal frame, in the class
\begin{equation*}
%\label{funct_class}
\mathcal{A} \,:=\, 
\Big\{(d_1,d_2,d_3)\,:\,(d_1|d_2|d_3)\in{\rm W}^{1,2}(I,{\rm SO}(3)),\, d_1^\prime\cdot d_2 = 0\mbox{ a.e. in } I
\Big\}.
\end{equation*}
In the above expression for $\mathscr E^{lim}$,
the energy density is defined as
\begin{equation}
\label{Qasmin}
\overline Q^{\,\theta}(\alpha,\beta):=\min_{\gamma\in\R}
\left\{c |M|^2 + 2c|\det M| + L^\theta(M)\,:\,M = \left( \begin{array}{rr} \alpha & \beta  \\ \beta  & \gamma  \end{array}\right)\right\}.
\end{equation}
We remark, overlooking for the moment the regularity 
of the deformations, that while in the 2D setting 
the admissible deformations are mappings $u$ from
the 2D domain $\omega$ to $\R^3$ such that
$u(\omega)$ is a developable surface, 
in the 1D limiting model they are curves $v$ 
endowed, at each point $v(x_1)$, with a 
orthonormal frame $(d_1(x_1),d_2(x_1),d_3(x_1))$, 
such that $v'=d_1$.
In particular, in the passage from 2D surfaces to 
``decorated'' curves, the isometry constraint 
($u(\omega)$ developable surface) is lost, 
but still recognizable in the constraint $d_1'\cdot d_2=0$
(the narrow strip cannot bend within its plane)
appearing in the definition of the admissible class
$\mathcal A$. 
Recall that the physical meaning of the 
relevant quantities 
$d_1'\cdot d_2$, $d_1^\prime\cdot d_3$,
and $d_2^\prime\cdot d_3$, appearing in the definition
of the limiting functional $\mathscr E^{lim}$,
is that of flexural strain around the thickness axis
(in short, \emph{flexure}), 
of flexural strain around the width axis,
and of torsional strain (in short, \emph{torsion}).  

\begin{theorem}
\label{cor:twist}
If $(u_\ep)\subset{\rm W}^{2,2}_{{\rm iso},\ep}(\omega,\R^3)$
is a sequence of minimisers of $\mathscr E^\ep$
then, up to a subsequence, there exist
$u\in{\rm W}^{2,2}(I,\R^3)$ and 
a minimiser $(d_1|d_2|d_3)\in\mathcal A$ of
$\mathscr E^{lim}$ 
with $d_1 = y^\prime$ 
such that 
\[
y_\ep\rightharpoonup y\ \mbox{ in }\ 
{\rm W}^{2,2}(\omega,\R^3),
\qquad
\nabla_\ep y_\ep \rightharpoonup (d_1|d_2)
\ \mbox{ in }\ {\rm W}^{1,2}(\omega,\R^{3\times2}),
\]
and
\begin{equation}
\label{here_gamma}
A_{y,\ep}\rightharpoonup \left( \begin{array}{rr} d_1^\prime\cdot d_3 & d_2^\prime\cdot d_3  \\ d_2^\prime\cdot d_3  & \gamma  \end{array}\right)\ \mbox{ in }\ 
{\rm L}^2(\omega,{\rm Sym}(2)),
\qquad\mbox{for some }\ \gamma\in\ {\rm L}^2(\omega,\R^3).
\end{equation}
Moreover,
\begin{equation*}
%\label{cvg_minima_T}
\min_{{\rm W}^{2,2}_{{\rm iso},\ep}(\omega,\R^3)}
\mathscr E^\ep
\longrightarrow
\min_\mathcal A
\mathscr E^{lim}.
\end{equation*}
\end{theorem}

The explicit expression of $\overline Q_T^{\,\theta}$ is 
\[
\overline Q_T^{\,\theta}(\alpha,\beta)
\,=\, 
\left\{\begin{array}{lc}
\frac{\mu}3 k\big(a_\theta\alpha - b_\theta\beta\big)
+ \frac{\mu}{12} k^2\Big(2-\frac{c_1}ca_\theta^2\Big) + \frac{\beta_T}2, 
& 
\quad\mbox{in}\quad\mathcal D_T\\
 \frac{\mu}3\big[(1+\gamma)\beta^2 - k b_\theta\beta\big]
+\frac{\mu}{12} k^2
\Big(2-\frac{a_\theta^2}{1+\gamma}\Big) + \frac{\beta_T}2,  
&
\quad\mbox{in}\quad\mathcal U_T\\
\frac{\mu}{12}(1+\gamma)
\frac{(\alpha^2+\beta^2)^2}{\alpha^2} + 
\frac{\mu}6 k
\Big(a_\theta\frac{\alpha^2-\beta^2}{\alpha} 
- 2b_\theta\beta\Big)
+\frac{\mu}6 k^2 + \frac{\beta_T}2,\quad  
& 
\quad\mbox{in}\quad\mathcal V_T, 
\end{array}\right.
\]
where $\gamma$ is given in \eqref{def:gamma},
$a_\theta$, $b_\theta$ and $\kk$ 
are defined as in \eqref{A^T},
and
\begin{align}
\mathcal D_T
&:=\left\{(\alpha,\beta)\in\R^2:
\frac{k\,a_\theta}{1+\gamma}
\,\alpha > \beta^2+\alpha^2\right\},\label{set_1}\\
\mathcal U_T
&:=\left\{(\alpha,\beta)\in\R^2:
\frac{k\,a_\theta}{1+\gamma}
\,\alpha \leq \beta^2-\alpha^2\right\},\label{set_2}\\
\mathcal V_T
&:=\R^2\setminus(\mathcal D_T\cup\mathcal U_T).\label{set_3}
\end{align}
Note that $\mathcal U_T$ is the interior of a 
(possibly degenerate) hyperbola. 
As for the set $\mathcal D_T$, note that 
it coincides with the interior of a disk,
whenever $\theta\notin\{\pi/4,3\pi/4\}$ (hence $a_\theta\neq0$),
while it reduces to the empty set
if $\theta=\pi/4$ or $\theta=3\pi/4$.
Note that $\overline Q_T^{\,\theta}$ is an affine
function in $\mathcal D_T$, and 
it is a parabola in $\beta$ in $\mathcal U_T$.
A deeper analysis of the 
above expression also shows
that $\overline Q_T^{\,\theta}$ is a continuous function.

\begin{figure}[htbp]
\begin{center}
\includegraphics[width=4.8cm]{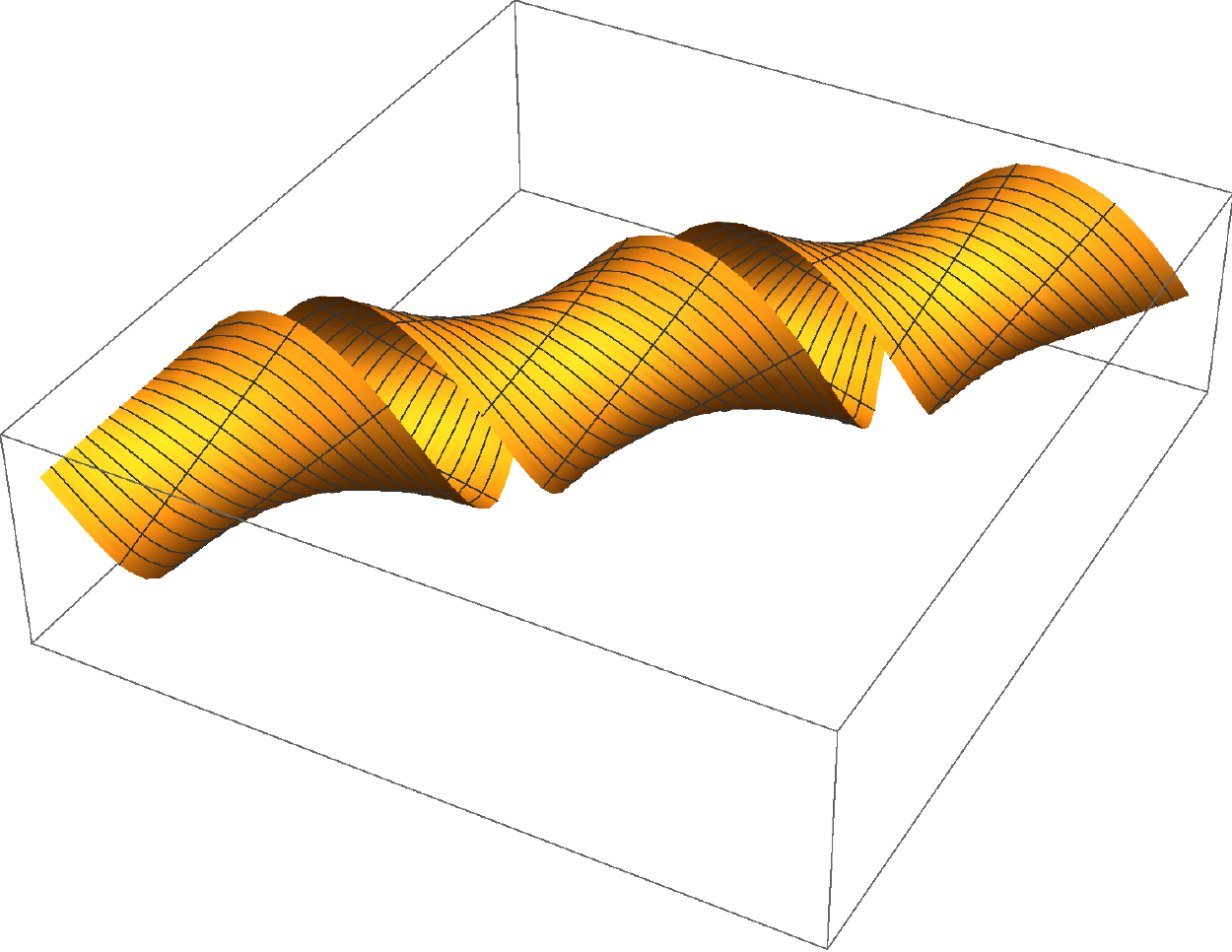}
\quad
\includegraphics[width=5.3cm]{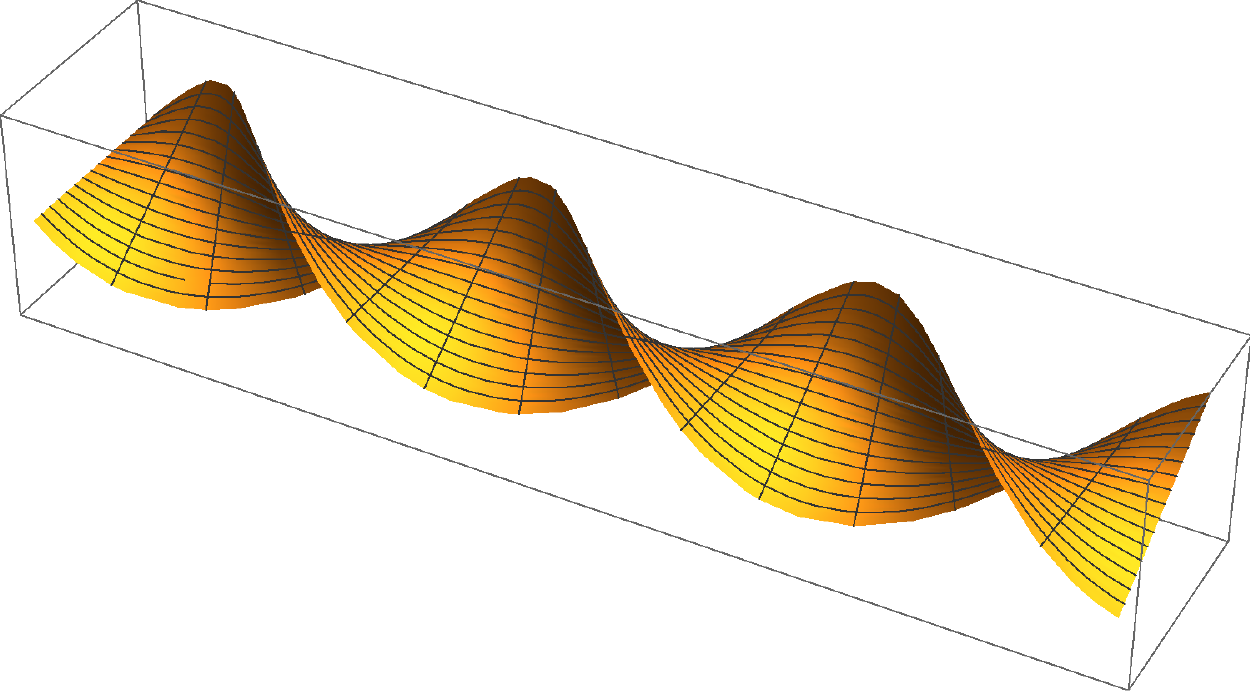}
\quad
\includegraphics[width=5.5cm]{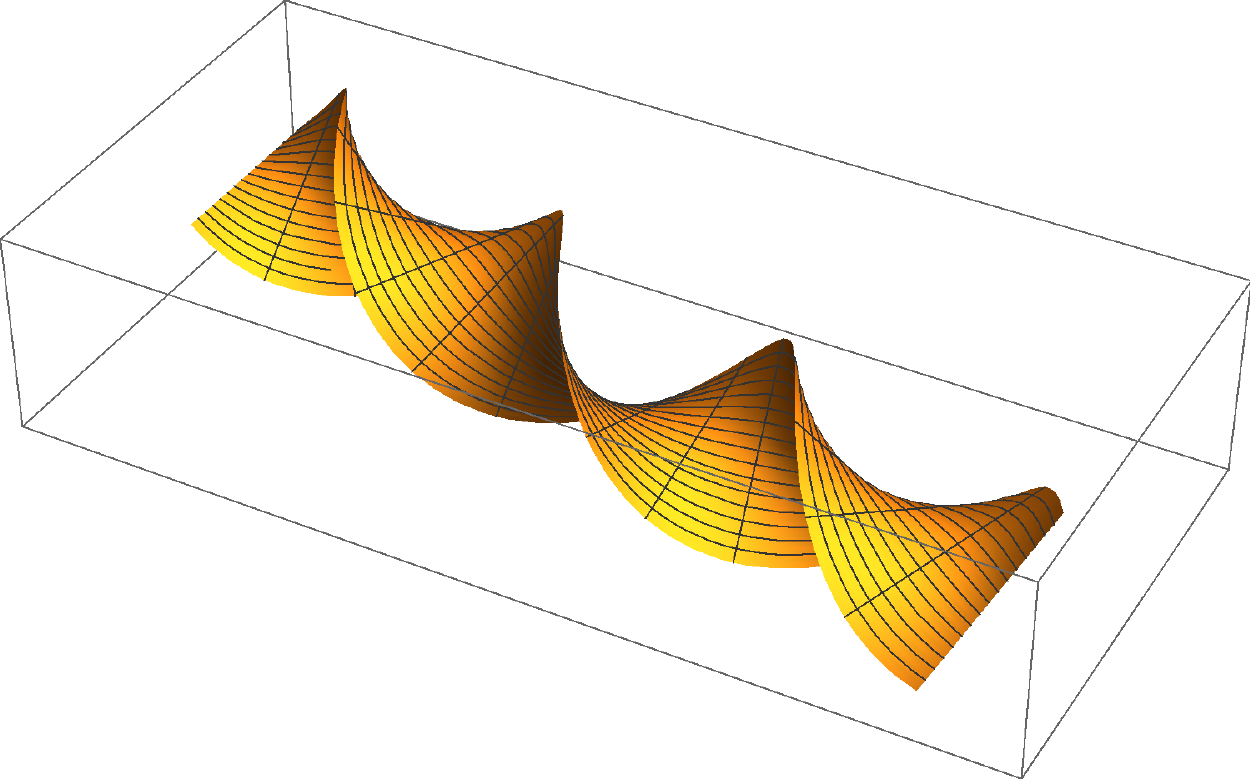}
\end{center}
\caption{
Minimal-energy configurations for the 1D rod model,
in the case $\theta=\pi/4$. The central figure
exhibits pure (constant) torsion (and zero flexure),
while the other plots are characterised by 
(constant) $\bar\alpha\neq0$ and $\bar\beta\neq0$.
In particular, in the rightmost picture, flexure
is close to zero. 
} 
\label{fig:config_min}
\end{figure}

Finally, one can show (see \cite[Lemma 3.3]{Ago_DeS_Kou})
that for every $0\leq\theta<\pi$,
$\overline Q_T^{\,\theta}$ attains its minimum value
precisely on the following subset of $\R^2$:
\begin{equation}
\label{min_set}
\left[-\frac{\kk}2
\left(\frac{1+\cos2\theta}{1+\gamma}\right),
\frac{\kk}2
\left(\frac{1-\cos2\theta}{1+\gamma}\right)\right]
\,\times\,
\left\{\frac{\kk}2
\left(\frac{\sin 2\theta}{1+\gamma}\right)\right\},
\qquad\qquad
\kk := \frac6{\pi^2}\frac{\alpha_0}{h_0}.
\end{equation}
Moreover,
\begin{equation*}
%\label{min_Q_T}
\min_{\R^2}\overline Q_T^{\,\theta}
\,=\,
\frac{\mu}{12}\kk^2\left(\frac{1+2\gamma}{1+\gamma}\right)
+\frac{\beta_T}2.
\end{equation*}
Building on this, 
we can construct minimum-energy configurations for 
the 1D model \eqref{funct_lim_T}--\eqref{Qasmin},
solving the system
\begin{equation*}
%\label{ODINA}
(d_1,d_2,d_3)\in\mathcal A,
\qquad\qquad
d_1^\prime\cdot d_3=\bar\alpha,
\qquad\qquad
d_2^\prime\cdot d_3=\bar\beta,  
\end{equation*}
where $\bar\alpha$ and $\bar\beta$ are two constant
values chosen in the set \eqref{min_set}.
Figure \ref{fig:config_min} shows a plot
of three minimal-energy configurations corresponding
to the case $\theta=\pi/4$. 
They both exhibit nontrivial flexure and torsion.

The previous analysis shows that there are many 
configurations realising the minimum of the 1D energy functional \eqref{funct_lim_T},
at fixed $\theta$.
It  would be interesting to try to derive more refined and ``selective'' 1D models capable of breaking this degeneracy, by discriminating between the configurations shown in Fig.~\ref{fig:config_min}.
We plan to address this issue in future work.  

%%%%%%%%%%%%%%%%%%%%%%%%%%%%%%%%%
%%%%%%%%%%%%%%%%%%%%%%%%%%%%%%%%%

\section{Some variants}
\label{sec:var}

%%%%%%%%%%%%%%%%%%%%%%%%%%%%%%%%%
%%%%%%%%%%%%%%%%%%%%%%%%%%%%%%%%%

Thin sheets of nematic elastomers with twist texture lead to a 2D model with negative target Gaussian curvature.
In this section,
we report on some variants of the 3D-to-2D derivation
of Subsections \ref{3D}--\ref{2D},  to explore some different scenarios. We show that nematic sheets with splay-bend texture lead to a 2D model with vanishing Gaussian target curvature (see Subsection \ref{SB}), while nematic sheets with director uniformly aligned, and perpendicular to the mid-surface, can lead to a 2D model with positive Gaussian target curvature (see Subsection \ref{Sharon}).
Finally, in Subsection \ref{bilayer}, we show that  thin bilayer sheets provide a rich model system in which positive, zero, or negative target Gaussian curvature can be produced at will, 
by tuning the spontaneous strains in the two halves of the bilayer.

%%%%%%%%%%%%%%%%%%%%%%%%%%%%%%%%%

\subsection{Splay-bend nematic elastomer sheets}
\label{SB}

%%%%%%%%%%%%%%%%%%%%%%%%%%%%%%%%%

\begin{figure}[htbp]
\begin{center}
\includegraphics[width=7cm]
{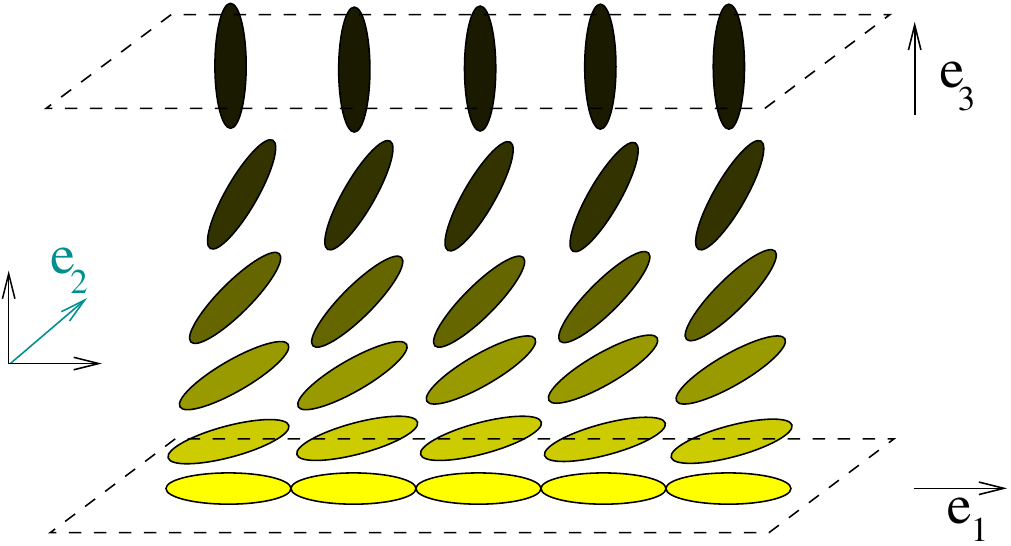}
\end{center}
\caption{
Sketch of the splay-bend director field.
} 
\label{fig:splay}
\end{figure}

In this subsection, we focus attention on
a \emph{splay-bend} nematic elastomer thin sheet.
In this case, the distribution of the nematic director 
along the thickness is given by
\begin{equation*}
\hat n_h(z_3)
\,:=\,
\left(
\begin{array}{c}
\cos\big(\frac{\pi}4+\frac\pi2\frac{z_3}h\big)\\
0\\
\sin\big(\frac{\pi}4+\frac\pi2\frac{z_3}h\big)
\end{array}
\right),
\qquad\qquad
z_3\in\Big[-\frac h2,\frac h2\,\Big],
\end{equation*}
see Figure \ref{fig:splay}.
Note that this distribution, which is constant
on each horizontal plane, coincides
with $\mathsf e_1$ on the bottom face
and with $\mathsf e_3$ on the top face.
Using again expression \eqref{prot_energy}
for the energy density,
the 3d-to-2d derivation is similar
%analogous,  mutatis mutandis, 
to that performed in
Subsections \ref{2D}--\ref{3D}. 
In particular, we arrive to (rescaled) spontaneous
strains of the form \eqref{eq:resc_SS},
where now
\begin{equation*}
B(x_3)
\,=\,
\frac12\frac{\alpha_0}{h_0}
\left(
\begin{array}{ccc}
\frac13-\cos^2\big(\frac{\pi}4+\frac\pi2x_3\big) 
& 0 & \frac12\sin\big(\frac{\pi}2+\pi x_3\big) \\
0 & \frac13 & 0\\
\frac12\sin\big(\frac{\pi}2+\pi x_3\big) 
& 0 &  \frac13- \sin^2\big(\frac{\pi}4+\frac\pi2x_3\big)\\
\end{array}
\right).
\end{equation*}
Recall that only the top left $2{\times}2$ part $\check B$
of $B$ enters into the definition of the
2D energy density $\overline Q_2$, see \eqref{barQ2}.
Considering that the relaxed density $Q_2$ has the
same expression as in \eqref{Q2}, 
after some computations 
(see \cite{Ag_De_bend} for more details) 
we obtain, similarly to the twist case, 
\[
\overline Q_2(G)
\,=\,
\alpha_{SB}\,Q_2(G-\bar A_{SB})+\beta_{SB},
\qquad\quad
\mbox{for every }\ G\in\R^{2\times 2}, 
\]
and in turn the limiting energy functional
\begin{equation*}
\mathscr F_{lim}^\ep(y)
\,=\,
\frac12
\int\limits_{\omega^\ep}
\left[
\alpha_{SB}\,Q_2(A_y(x')-\bar A_{SB})+\beta_{SB}
\right]
{\rm d}x',
\qquad\quad y\in{\rm W}^{2,2}_{\rm iso}(\omega^\ep,\R^3),
\end{equation*}
where
\begin{equation*}
%\label{const_SB}
\alpha_{SB}\,=\,
\alpha_T
\,=\,\frac1{12},
\quad\qquad
\bar A_{SB}\,=\,\frac{6}{\pi^2}\frac{\alpha_0}{h_0}
\,{\rm diag(-1,0)},
\quad\qquad
\beta_{SB}\,=\,\mu\,(1+\gamma)
\Big(\frac{\pi^4-12}{16}\Big)
\frac{\alpha_0^2}{h_0^2}.
\end{equation*}
%and $f_{lim}^\ep$ is defined as in \eqref{f_lim}.
Moreover, $\gamma$ is given in terms of the 3D
parameters by \eqref{def:gamma}.
We refer the reader to \cite{Ago_DeS_Kou} for the
derivation of a corresponding 1D model as $\ep\downarrow0$,
in the spirit of Subsection \ref{1D}.

If $\{v_h\}$ a low-energy sequence of (physical) deformations,
with corresponding energies 
$\widehat{\mathscr F}_h^\ep(v_h)$,
we have that
\begin{multline*}
\widehat{\mathscr F}_{h_0}^\ep(v_{h_0})
\,\cong\,
\min_{y\in{\rm W}^{2,2}_{\rm iso}(\omega^\ep,\R^3)}
\frac{\mu\,h_0^3}{12}\int\limits_{\omega^\ep}
\!\!\left\{\bigg|{\rm A}_y(x')\!-\!
\frac{6}{\pi^2}\frac{\alpha_0}{h_0}\,
{\rm diag}\big(-1,0\big)\bigg|^2
\!\!+\gamma
\Big({\rm H}_y(x')\!+\!
\frac{6}{\pi^2}\frac{\alpha_0}{h_0}\Big)^2\!\right\}\!\dd x'\\ 
+\, \mu\,(1+\gamma)\,h_0\,\alpha_0^2
\Big(\frac{\pi^4-12}{32}\Big)\!|\omega^\ep|,
\end{multline*}
for a given small thickness $h_0$,
where the approximate identity holds 
modulo terms of order higher than $3$ in $h_0$.

Observe that
$\min\mathscr F_{lim}^\ep=\beta_{SB}|\omega^\ep|/2$, since,
differently from the twist case (cf.\ Remark \ref{frustration}),
the function $x'\mapsto Q_2(A_y(x')-\bar A_{SB})$
can be minimised to the constant value zero,
by any $y\in{\rm W}^{2,2}_{\rm iso}(\omega^\ep,\R^3)$
such that $A_y\equiv\bar A_{SB}$.

%%%%%%%%%%%%%%%%%%%%%%%%%%%%%%%%%

\subsection{Constant director along the thickness}
\label{Sharon}

%%%%%%%%%%%%%%%%%%%%%%%%%%%%%%%%%

In this subsection, we consider the case
where the nematic director is constant along the thickness
and the dependence of the 
spontaneous strain \eqref{nem_tens} on the 
thickness variable is encoded by the 
magnitude parameter $a$.
More precisely,
using the notation of Subsection \ref{3D},
we suppose that the nematic director
$\hat n(z_3)$ is constantly equal to some
$\mathsf n\in\Sph^2$, whereas the (constant) parameter
$a_h$ in \eqref{spont_phys} is here replaced
by the function
\begin{equation*}
\hat a_h(z_3)
\,:=\,
1+\frac{\alpha_0}{h_0}z_3,
\qquad\qquad
z_3\in\Big(\!-\frac h2,\frac h2\Big).
\end{equation*}
All in all, the (physical) spontaneous strain
of this system is defined as
\begin{equation*}
\widehat C_h(z_3)
\,:=\,
\hat a_h^{2/3}(z_3)\mathsf n\otimes\mathsf n
+\hat a_h^{-1/3}(z_3)\big({\rm I}-\mathsf n\otimes\mathsf n\big).
\end{equation*}
Modelling the system using again the prototypical
energy density \eqref{prot_energy} and following
the same notation and steps of Subsections \ref{3D}--\ref{2D}
we arrive at the rescaled energy densities $W_h(x_3,\cdot)$,
$x_3\in(-1/2,1/2)$,
characterised by the (rescaled) spontaneous strains
\begin{align*}
\overline C_h(x_3)
\,:=\,
\widehat C_h(hx_3)
&
\,=\,\Big(1+\frac{\alpha_0}{h_0}z_3\Big)^{2/3}
\!\mathsf n\otimes\mathsf n
+
\Big(1+\frac{\alpha_0}{h_0}z_3\Big)^{-1/3}
\big({\rm I}-\mathsf n\otimes\mathsf n\big)\nonumber\\
&
\,=\,{\rm I}-2\,h\,B(x_3)+R^h(x_3),
\qquad\quad
B(x_3)
\,:=\,
\frac{x_3}2\frac{\alpha_0}{h_0}
\left(\frac{\rm I}3-\mathsf n\otimes\mathsf n\right),
\end{align*}
where $\|R^h\|_{\infty}=o(h)$.
The corresponding 2D model that we obtain in the end
associates with each 
$y\in{\rm W}^{2,2}_{\rm iso}(\omega^\ep,\R^3)$ the energy
\begin{equation}\label{001}
\mathscr F_{lim}^\ep(y)
\,=\,
\frac12
\int\limits_{\omega^\ep}
\overline Q_2(A_y(x'))
{\rm d}x'
\,=\,
\frac1{24}
\int\limits_{\omega^\ep}
Q_2\left(A_y(x')-\check M\right)
{\rm d}x',
\end{equation}
where the $2{\times}2$ symmetric matrix $\check M$
is given by the formula
\begin{equation*}
\check M
\,=\,
\frac12\frac{\alpha_0}{h_0}
\left[(\mathsf n\otimes\mathsf n)^{\check{}}-
\frac{{\rm I}_2}3\right],
\end{equation*}
and 
$(\mathsf n\otimes\mathsf n)^{\check{}}$
is the $2{\times}2$ upper left part of 
$\mathsf n\otimes\mathsf n$.
For example, in the case where  
$\mathsf n=\mathsf e_3$, we have that
\[
\check M
\,=\,
\left(
\begin{array}{ccc}
m_0 & 0\\
0 & m_0
\end{array}
\right),
\qquad\qquad
m_0\,:=\,
-\frac{\alpha_0}{6h_0}.
\]
Notice that this corresponds to a target curvature
with positive Gaussian curvature.
However, observable minimal-energy configurations,
will always exhibit zero Gaussian curvature,
because of the isometry constraint
they are subjected to. The reader is referred to \cite{Ag_De_bend} for more details, and for plots of minimal-energy configurations predicted by \eqref{001}.

%%%%%%%%%%%%%%%%%%%%%%%%%%%%%%%%%

\subsection{Bilayers}
\label{bilayer}

%%%%%%%%%%%%%%%%%%%%%%%%%%%%%%%%%

Finally, we consider the case of a bilayer governed, again,
by the prototypical energy density \eqref{nostra_W_0}.
A similar problem has been considered in 
\cite{Schmidt2007} and also in \cite{Bartels}.
More precisely, we consider a model for bilayers
where in the physical reference configuration
$\Omega_h^\ep$ the energy density 
$\widehat W_h=\widehat W_h(z_3,F)$
is given by
\eqref{W_h_W_0}--\eqref{nostra_W_0} and
the spontaneous (right Cauchy-Green) strain is defined as 
\begin{equation*}
\widehat C_h(z_3)
\,:=\,
{\rm I}\,-\,2\,h\,\widehat B(z_3),
\qquad\qquad
\widehat B(z_3)
\,:=\,
\left\{
\begin{array}{rl}
M_1 & \mbox{\ if }z_3\in[0,h/2)\\
M_1 & \mbox{\ if }z_3\in
(-h/2,0),
\end{array}
\right.
\end{equation*}
for some fixed $M_1,M_2\in{\rm Sym}(3)$,
whose physical dimension is that of inverse length.
We again use the notation
$\widehat{\mathscr F}_h^\ep$ for 
the 3D (physical) free-energy
functional, defined as in \eqref{phys_quan}.
Passing to the rescaled reference configuration 
and the corresponding  rescaled energy densities
$W_h:(-1/2,1/2)\times\R^{3\times3}\longrightarrow[0,+\infty]$
defined as $W_h(x_3,F):=\widehat W_h(hx_3,F)$,
we have that
$W_h(x_3,F)=W_0\big(F\overline U_h(x_3)^{-1}\big)$,
where $W_0$ is defined as in \eqref{nostra_W_0} and
$\overline U_h(x_3):=\widehat C_h(hx_3)^{1/2}$,
so that
\begin{equation*}
\overline U_h(x_3)^{-1}
\,=\,
{\rm I}+h\,B(x_3)+o(h),
\qquad\mbox{with}\quad
B(x_3)
\,:=\,
\widehat B(h\,x_3)
\,=\,
\left\{
\begin{array}{rl}
M_1 & \mbox{\ if }x_3\in[0,1/2)\\
M_2 & \mbox{\ if }x_3\in
(-1/2,0).
\end{array}
\right.
\end{equation*}
The expression above shows that the spontaneous strains are compatible with the format used in \cite{Schmidt2007}, so that we can proceed as in Section 2 for the derivation of 
a limiting 2D plate model. The model is nonlinear, in the sense that it covers the regime of arbitrarily large bending deformations. The associated strains are however always small, in view of the thinness of the sheet and the fact that, in the thin limit, the mid-plane deforms isometrically. In particular, the spontaneous strains are small. Since $\overline U_h= I-hB + o(h)$, we have that the spontaneous linear (or infinitesimal) strains in the bilayer are described by
\begin{equation}\label{ads}
\widehat E_h(z_3)
\,=\,
-\,h\,\widehat B(z_3)
\,=\,
\left\{
\begin{array}{rl}
E_1 := \,- \,h M_1 & \mbox{\ if }z_3\in[0,h/2)\\
E_2 := \,- \,h M_2 & \mbox{\ if }z_3\in
(-h/2,0),
\end{array}
\right.
\end{equation}
where $E_1$ and $E_2$ are the spontaneous linear strains in the top and bottom half of the bilayer, respectively.

Proceeding as in Subsection \ref{2D}, to which we refer for the notation, 
we can use again the results of \cite{Schmidt2007}
to deduce from the functional 
$h^{-2}\int_{\Om^\ep}W_h(x_3,\na_h(y)){\rm d}x_3$,
via compactness and $\Gamma$-convergence arguments,
the 2D limit functional
\begin{equation}
\mathscr F_{lim}^\ep(y)
\,:=\,
\frac12
\int\limits_{\omega^\ep}
\overline Q_2\big(A_y(x')\big)
{\rm d}x',
\qquad
y\in{\rm W}^{2,2}_{\rm iso}(\omega^\ep,\R^3),
\end{equation}
where
\[
\overline Q_2
\,:=\,
\min_{D\in\R^{2\times2}}\int\limits_{-1/2}^{1/2}
Q_2\big(D+t\,G+\check B(t)\big){\rm d}t,
\] 
and with $Q_2$ and $\check B$ defined as in 
\eqref{Q2}--\eqref{def:gamma} and \eqref{cek_B},
respectively.
Observe that, denoting by $L$ the bilinear
form associated with $Q_2$, namely, 
\[
L(G,H)\,:=\,
2\mu\big(\sym G\cdot\sym H+\gamma\,\tr G\,\tr H\big),
\qquad\quad
G,H\in\R^{2\times2},
\]
we have that
\begin{align}
\overline Q_2(G)
&\,=\,
\min_{D\in\R^{2\times2}}\int\limits_{-1/2}^{1/2}
\Big[Q_2(D)+Q_2\big(t\,G+\check B(t)\big)
+2L\big(D,t\,G+\check B(t)\big)\Big]{\rm d}t,
\nonumber\\
&\,=\,
\int\limits_{-1/2}^{1/2}
Q_2\big(t\,G+\check B(t)\big)\,{\rm d}t
+
\min_{D\in\R^{2\times2}}
\left[Q_2(D)
+2L\left(D,\int_{-\frac12}^{\frac12}\big(t\,G+\check B(t)\big){\rm d}t\right)\right].
\label{Q_2+min}
\end{align}
Now, simple computations give that
\begin{align}
\int\limits_{-1/2}^{1/2}
Q_2\big(t\,G+\check B(t)\big)\,{\rm d}t
&\,=\,
Q_2(G)\int\limits_{-1/2}^{1/2}t^2\,{\rm d}t
+\int\limits_{-1/2}^{1/2}
Q_2\big(\check B(t)\big)\,{\rm d}t
+
2L\left(G,\int\limits_{-1/2}^{1/2}t\check B(t)
\,{\rm d}t\right),\nonumber\\
&\,=\,
\frac1{12}
Q_2(G)+
\frac14L(G,\check M_1-\check M_2)+
\frac12\big[Q_2(\check M_1)+Q_2(\check M_2)\big],\nonumber\\
&\,=\,
\frac1{12}Q_2\Big(G+\frac32(\check M_1-\check M_2)\Big)+
c,\label{etc+c}
\end{align}
where in the second equality we have used the fact that
\[
\check B(t)
\,=\,
\left\{
\begin{array}{rl}
\check M_1 & \mbox{\ if }t\in[0,1/2)\\
\check M_2 & \mbox{\ if }t\in
(-1/2,0),
\end{array}
\right.
\]
and where we have set
\[
c\,:=\,
\frac1{16}\Big[
5\,Q_2(\check M_1)
+
5\,Q_2(\check M_2)
+
6\,L(\check M_1,\check M_2)
\Big].
\]
Also, note that the minimum problem
in \eqref{Q_2+min} does not depend on $G$,
since $\int_{-1/2}^{1/2}t{\rm d}t=0$,
so that it is equivalent to
$\min_{D\in\R^{2\times 2}}f^{\check B}(D)$,
where 
\begin{equation}
\label{def_f}
f^{\check B}(D)
\,:=\,
Q_2(D)
+2L(D,\check N),
\qquad\qquad
\check N\,:=\,
\int\limits_{-\frac12}^{\frac12}\check B(t){\rm d}t.
\end{equation}
It is easy to show that
\begin{equation}
\label{min_f}
\min_{D\in\R^{2\times 2}}f^{\check B}(D)
\,=\,
f^{\check B}(-\check N)
\,=\,
-Q_2(\check N)
\,=\,
-\frac14 Q_2(\check M_1+\check M_1),
\end{equation}
where in the last equality we have used the
explicit expression of $\check B$, 
and in turn that of $\check N$, 
in terms of the matrices $\check M_i$'s.
Putting together \eqref{Q_2+min},
\eqref{etc+c}, and \eqref{def_f}--\eqref{min_f},
we obtain that
\[
\overline Q_2(G)
\,=\,
\frac1{12}
Q_2\Big(G+\frac32(\check M_1-\check M_2)\Big)
+c
-\frac14 Q_2(\check M_1+\check M_1)
\,=\,
\frac1{12}
Q_2\Big(G+\frac32(\check M_1-\check M_2)\Big)
-\frac1{16}
Q_2(\check M_1+\check M_2).
\]
All in all, turning back to the physical variables
and considering a low-energy sequence of deformations $\{v_h\}$,
we have that the limiting 2D plate theory 
can be summarised by the approximate identity 
\begin{align*}
\widehat{\mathscr F}_{h_0}^\ep(v_{h_0})
&\,\cong\,
\min_{y\in{\rm W}^{2,2}_{\rm iso}(\omega^\ep,\R^3)}
\frac{h_0^3}{2}\int\limits_{\omega^\ep}
\overline Q_2\big({\rm A}_y(x')\big)\dd x'\\
&\,=\,
\min_{y\in{\rm W}^{2,2}_{\rm iso}(\omega^\ep,\R^3)}
\frac{h_0^3}{24}
\int\limits_{\omega^\ep}
Q_2\Big({\rm A}_y(x')-\frac{3}{2h_0}(\check E_1-\check E_2)\Big)
\dd x'
-
\frac{h_0|\omega^\ep|}{32}
Q_2(\check E_1+\check E_2),
\end{align*}
where $Q_2$ is defined as in 
\eqref{Q2}--\eqref{def:gamma}, and  $\check E_1$  and $\check E_2$ denote the sub-matrices consisting of the first two rows and columns of the spontaneous linear strains in the top and bottom part of the bilayer, respectively, see \eqref{ads}. Depending on the values the components of $\check E_1$  and $\check E_2$, one can obtain target curvature tensors $3(\check E_1-\check E_2)/2h_0$
with positive, zero, or negative determinant, and hence positive, zero, or negative target Gaussian curvature. In all cases, however, the isometry constraint will have the effect that energy minimising configurations will be ribbons wrapped around a cylinder.

%%%%%%%%%%%%%%%%%%%%%%%%%%%%%%%%%%%%%%%%%%

\section{Conclusions and discussion}

%%%%%%%%%%%%%%%%%%%%%%%%%%%%%%%%%%%%%%%%%%

In this paper, we have discussed a strategy to derive
a model for nematic elastomer thin and narrow films,
with a twist nematic texture.
The strategy uses a rigorous dimension reduction argument,
from 3D to 2D (plate model),
and then from 2D to 1D (narrow ribbon model).
This procedure leads to a degenerate model,
which admits as minimisers both spiral ribbons 
and helicoid-like shapes.
It would be interesting to derive more selective
1D models, able to discriminate
between these two types of configurations, and to deliver unique minimisers in different regimes of the relevant material and geometric parameters governing the behaviour of thin and narrow strips. 
This will be attempted in future work.
We refer the reader to \cite{Tomassetti} 
for further discussion of spiral ribbon vs. 
helicoid-like shapes in thin sheets of nematic elastomers.

In addition, we have considered variants to thin nematic sheets with twist texture, which leads to 2D models with negative Gaussian target curvature.
In the context of nematic elastomers, these are the splay-bend texture and uniform alignment of the director perpendicular to the mid-plane, which lead to zero and positive Gaussian target curvature, respectively. We have also shown that thin bilayer sheets provide a rich model system in which positive, zero, or negative target Gaussian curvature can be produced at will, 
by tuning the spontaneous strains in the two halves of the bilayer.

%%%%%%%%%%%%%%%%%%%%%%%%%%%%%%%%%%%%%%%%%%%%%%%
%%%%%%%%%%%%%%%%%%%%%%%%%%%%%%%%%%%%%%%%%%%%%%%

\subsection*{Acknowledgements} We thank E. Sharon for valuable discussions.
We gratefully acknowledge the support by the European Research Council through the ERC Advanced Grant 340685-MicroMotility.

%%%%%%%%%%%%%%%%%%%%%%%%%%%%%%%%%%%%%%%%%%%%%%%
%%%%%%%%%%%%%%%%%%%%%%%%%%%%%%%%%%%%%%%%%%%%%%%

\medskip
\noindent
The authors declare that they have no conflict of interest.

%%%%%%%%%%%%%%%%%%%%%%%%%%%%%%%%%%%%%%%%%%%%%%%%
%%%%%%%%%%%%%%%%%%%%%%%%%%%%%%%%%%%%%%%%%%%%%%%%

\bibliographystyle{plain}

%%%%%%%%%%%%%%%%%%%%%%%%%%%%%%%%%%%%%%%%%%%%%%%%
%%%%%%%%%%%%%%%%%%%%%%%%%%%%%%%%%%%%%%%%%%%%%%%%

\end{document}